\input amstex
\input amsppt.sty
\magnification=\magstep1
\hsize=33.5truecc
\vsize=23truecm
\baselineskip=16truept
\NoBlackBoxes
\TagsOnRight \pageno=1 \nologo
\def\Z{\Bbb Z}
\def\N{\Bbb N}

\def\l{\left}
\def\r{\right}
\def\bg{\bigg}
\def\({\bg(}
\def\[{\bg\lfloor}
\def\){\bg)}
\def\]{\bg\rfloor}
\def\t{\text}
\def\f{\frac}

\def\bi{\binom}
\def\eq{\equiv}

\def\ls{\leqslant}
\def\gs{\geqslant}
\def\mo{\roman{mod}}

\def\al{\alpha}

\def\Proof{\noindent{\it Proof}}

\def\Remark{\medskip\noindent{\it  Remark}}

\def\Ack{\medskip\noindent {\bf Acknowledgment}}
\hbox {J. Number Theory 183 (2018), 146--171.}
\bigskip
\topmatter
\title Arithmetic properties of Delannoy numbers and Schr\"oder numbers\endtitle
\author Zhi-Wei Sun\endauthor
\leftheadtext{Zhi-Wei Sun}
\rightheadtext{Properties of Delannoy numbers and Schr\"oder numbers}
\affil Department of Mathematics, Nanjing University\\
 Nanjing 210093, People's Republic of China
  \\  zwsun\@nju.edu.cn
  \\ {\tt http://maths.nju.edu.cn/$\sim$zwsun}
\endaffil
\abstract Define
$$D_n(x)=\sum_{k=0}^n\binom nk^2x^k(x+1)^{n-k}\ \ \ \t{for}\ n=0,1,2,\ldots$$
and
$$s_n(x)=\sum_{k=1}^n\frac1n\binom nk\binom n{k-1}x^{k-1}(x+1)^{n-k}\ \ \ \t{for}\ n=1,2,3,\ldots.$$
Then $D_n(1)$ is the $n$-th central Delannoy number $D_n$,
and $s_n(1)$ is the $n$-th little Schr\"oder number $s_n$.
In this paper we obtain some surprising arithmetic properties of $D_n(x)$ and $s_n(x)$. We show that
$$\f1n\sum_{k=0}^{n-1}D_k(x)s_{k+1}(x)\in\Z[x(x+1)]\ \quad\t{for all}\ n=1,2,3,\ldots.$$
Moreover, for any odd prime $p$ and $p$-adic integer $x\not\equiv0,-1\pmod p$, we establish the supercongruence
$$\sum_{k=0}^{p-1}D_k(x)s_{k+1}(x)\equiv0\pmod{p^2}.$$
As an application we confirm Conjecture 5.5 in [S14a], in particular we prove that
$$\frac1n\sum_{k=0}^{n-1}T_kM_k(-3)^{n-1-k}\in\Z\quad\t{for all}\ n=1,2,3,\ldots,$$
where $T_k$ is the $k$-th central trinomial coefficient and $M_k$ is the $k$-th Motzkin number.
\endabstract
\thanks 2010 {\it Mathematics Subject Classification}. \,Primary 11A07, 11B75;
Secondary  05A10, 05A19, 11B65, 11C08, 33F99.
\newline\indent {\it Keywords}. Congruence, Polynomial, Delannoy number, Schr\"oder number, Motzkin number
\newline\indent Supported by the National Natural Science
Foundation of China (grant 11571162).
\endthanks

\endtopmatter
\document

\heading{1. Introduction}\endheading

For $m,n\in\N=\{0,1,2,\ldots\}$, the Delannoy number
$$D_{m,n}:=\sum_{k\in\N}\bi mk\bi nk2^k\tag1.1$$
in combinatorics counts lattice paths from $(0,0)$ to $(m,n)$ in which only east $(1, 0)$, north $(0, 1)$, and northeast $(1, 1)$ steps are allowed
(cf. R. P. Stanley [St99, p. 185]). The $n$-th central Delannoy number $D_n=D_{n,n}$ has another well-known expression:
$$D_n=\sum_{k=0}^n\bi nk\bi{n+k}k=\sum_{k=0}^n\bi{n+k}{2k}\bi{2k}k.\tag1.2$$

For $n\in\Z^+=\{1,2,3,\ldots\}$, the $n$-th little Schr\"oder number is given by
$$s_n:=\sum_{k=1}^nN(n,k)2^{k-1}\tag1.3$$
with the Narayana number $N(n,k)$ defined by
$$N(n,k):=\f1n\bi nk\bi n{k-1}\in\Z.$$
(See [Gr, pp.\,268--281] for certain combinatorial interpretations of the Narayana number $N(n,k)=N(n,n+1-k)$.)
For $n\in\N$, the $n$-th large Schr\"oder number is given by
$$S_n:=\sum_{k=0}^n\bi nk\bi{n+k}k\f1{k+1}=\sum_{k=0}^n\bi{n+k}{2k}C_k,\tag1.4$$
where $C_k$ denotes the Catalan number $\bi{2k}k/(k+1)=\bi{2k}k-\bi{2k}{k+1}$.
It is well known that $S_n=2s_n$ for every $n=1,2,3,\ldots$. Both little Schr\"oder numbers and large Schr\"oder numbers
have many combinatorial interpretations (cf. [St97] and [St99, pp.\,178, 239-240]); for example, $s_n$
is the number of ways to insert parentheses into an expression of $n+1$ terms with two or more items within a parenthesis, and $S_n$ is the
number of lattice paths from the point $(0,0)$ to $(n,n)$ with steps
$(1,0),(0,1)$ and $(1,1)$ which never rise above the line $y=x$.

Surprisingly, the central Delannoy numbers and Schr\"oder numbers arising naturally in enumerative combinatorics,
have nice arithmetic properties. In 2011 the author [S11b] showed that
$$\sum_{k=1}^{p-1}\f{D_k}{k^2}\eq(-1)^{(p-1)/2}2E_{p-3}\pmod p\ \ \t{and}\ \ \sum_{k=1}^{p-1}\f{S_k}{6^k}\eq0\pmod p$$
for any prime $p>3$, where $E_0,E_1,E_2,\ldots$ are the Euler numbers. In 2014 the author [S14a] proved that
$$\f1{n^2}\sum_{k=0}^{n-1}(2k+1)D_k^2\in\Z\quad\t{for all}\ n=1,2,3,\ldots,$$
and that
$$\sum_{k=0}^{p-1}D_k^2\eq\l(\f 2p\r)\pmod{p}\quad\t{for any odd prime}\ p,$$
where $(\f{\cdot}p)$ denotes the Legendre symbol.
\medskip

{\it Definition} 1.1. We define
$$D_n(x):=\sum_{k=0}^n\bi nk^2x^k(x+1)^{n-k}\ \quad\t{for}\ n\in\N,\tag1.5$$
and
$$s_n(x):=\sum_{k=1}^nN(n,k)x^{k-1}(x+1)^{n-k}\ \quad\t{for}\ n\in\Z^+.\tag1.6$$

Obviously $D_n(1)=D_n$ for $n\in\N$, and $s_n(1)=s_n$ for $n\in\Z^+$.

In this paper we obtain somewhat curious results involving the polynomials $D_n(x)$ and $s_n(x)$. Our first theorem is as follows.

\proclaim{Theorem 1.1} {\rm (i)} For any $n\in\Z^+$, we have
$$\f1n\sum_{k=0}^{n-1}D_k(x)s_{k+1}(x)=W_n(x(x+1)),\tag1.7$$
where
$$W_n(x)=\sum_{k=1}^nw(n,k)C_{k-1}x^{k-1}\in\Z[x]\tag1.8$$
with
$$w(n,k)=\f1k\bi{n-1}{k-1}\bi{n+k}{k-1}\in\Z.\tag1.9$$

{\rm (ii)} Let $p$ be an odd prime. For any $p$-adic integer $x$, we have
$$\aligned&\sum_{k=0}^{p-1}D_k(x)s_{k+1}(x)
\\\eq&\cases p(1-x(x+1))\pmod{p^3}&\t{if}\ x\eq0,-1\pmod p,
\\2p^2+\f{2x+1}{x(x+1)}p^2(x^2q_p(x)-(x+1)^2q_p(x+1))\pmod{p^3}&\t{otherwise},\endcases
\endaligned\tag1.10$$
where $q_p(z)$ denotes the Fermat quotient $(z^{p-1}-1)/p$
for any $p$-adic integer $z\not\eq0\pmod p$.
\endproclaim
\Remark\ 1.1. It is interesting to compare the new numbers $w(n,k)$ with the Narayana numbers $N(n,k)$.
\medskip

Clearly, Theorem 1.1 in the case $x=1$ yields the following consequence.

\proclaim{Corollary 1.1} For any positive integer $n$, we have
$$\f1n\sum_{k=0}^{n-1}D_ks_{k+1}=\sum_{k=1}^nw(n,k)C_{k-1}2^{k-1}\eq1\pmod 2.\tag1.11$$
Also, for any odd prime $p$ we have
$$\sum_{k=0}^{p-1}D_ks_{k+1}\eq2p^2(1-3q_p(2))\pmod{p^3}.\tag1.12$$
\endproclaim
\Remark\ 1.2. For the prime $p=588811$, we have $q_p(2)\eq1/3\pmod p$ and hence $\sum_{k=0}^{p-1}D_ks_{k+1}\eq0\pmod{p^3}$.
In 2016 J.-C. Liu [L] confirmed the author's conjecture (cf. [S11b, Conjecture 1.1]) that
$$\sum_{k=1}^{p-1}D_kS_k\eq-2p\sum_{k=1}^{p-1}\f{(-1)^k+3}k\pmod{p^4}\quad\t{for any prime}\ p>3.$$
\medskip

From Theorem 1.1 we can deduce a novel combinatorial identity.
\proclaim{Corollary 1.2} For any $n\in\Z^+$, we have
$$\sum_{k=1}^n\bi nk\bi{n+k}{k-1}\f{C_{k-1}}{(-4)^{k-1}}=\f{\lfloor(n+1)/2\rfloor}{4^{n-1}}\bi n{\lfloor n/2\rfloor}^2.\tag1.13$$
\endproclaim
\Remark\ 1.3. If we let $a_n$ denote the left-hand side of (1.13), then the Zeilberger algorithm  (cf. [PWZ, pp.\,101-119]) cannot find a closed form for $a_n$,
and it only yields the following second-order recurrence relation:
$$(n+1)^2a_n+(2n+3)a_{n+1}-(n+1)(n+3)a_{n+2}=0\ \ \t{for}\ n=1,2,3,\ldots.$$
\smallskip

Now we give our second theorem which can be viewed as a supplement to Theorem 1.1.
\proclaim{Theorem 1.2} Let $p$ be any odd prime.
Then
$$\sum_{k=0}^{p-1}kD_k(x)s_{k+1}(x)\eq 2(x(x+1))^{(p-1)/2}\pmod p.\tag1.14$$
In particular,
$$\sum_{k=0}^{p-1}kD_ks_{k+1}\eq2\l(\f 2p\r)\pmod p.\tag1.15$$
\endproclaim

In the next section we are going to show Theorems 1.1-1.2 and Corollary 1.2.
In Section 3 we will give applications of Theorems 1.1-1.2 to central trinomial coefficients, Motzkin numbers, and their generalizations.
Section 4 contains two related conjectures.

Throughout this paper, for any polynomial $P(x)$ and $n\in\N$, we use $[x^n]P(x)$ to denote the coefficient of $x^n$ in $P(x)$.

\heading{2. Proofs of Theorems 1.1-1.2 and Corollary 1.2}\endheading

Recall the following definition given in [S12a] motivated by the large Schr\"oder numbers.
\medskip

{\it Definition} 2.1. For $n\in\N$ we set
$$S_n(x)=\sum_{k=0}^n\bi nk\bi{n+k}k\f{x^k}{k+1}=\sum_{k=0}^n\bi{n+k}{2k}C_kx^k.\tag2.1$$

\proclaim{Lemma 2.1} We have
$$D_n(x)=\sum_{k=0}^n\bi nk\bi{n+k}kx^k\quad\t{for}\ n\in\N.\tag2.2$$
Also, for any $n\in\Z^+$ we have
$$D_{n+1}(x)-D_{n-1}(x)=2x(2n+1)S_n(x)\tag2.3$$
and
$$(x+1)s_n(x)=S_n(x).\tag2.4$$
\endproclaim
\Proof. For $k,n\in\N$, we obviously have
$$\align [x^k]D_n(x)=&[x^k]\sum_{j=0}^n\bi nj^2x^j(x+1)^{n-j}=\sum_{j=0}^k\bi nj^2\bi{n-j}{k-j}
\\=&\bi nk\sum_{j=0}^k\bi nj\bi k{k-j}=\bi nk\bi{n+k}k
\endalign$$
with the help of the Chu-Vandermonde identity (cf. [G, (3.1)]). This proves (2.2).

Now fix $n\in\Z^+$.  For $k\in\N$, by (2.2) we clearly have
$$\align&[x^{k+1}](D_{n+1}(x)-D_{n-1}(x))
\\=&\bi{n+1+(k+1)}{2(k+1)}\bi{2(k+1)}{k+1}-\bi{n-1+(k+1)}{2(k+1)}\bi{2(k+1)}{k+1}
\\=&\f{2n+1}{2k+1}\bi{n+k}{2k}\bi{2k+2}{k+1}=\f{2(2n+1)}{k+1}\bi{n+k}{2k}\bi{2k}k
\\=&[x^{k+1}]2x(2n+1)S_n(x).
\endalign$$
So (2.3) follows.

For each $k\in\N$, it is apparent that
$$\align [x^k](x+1)s_n(x)=&[x^k](x+1)\sum_{j=1}^nN(n,j)x^{j-1}(x+1)^{n-j}
\\=&\sum_{0<j\ls k}N(n,j)\bi{n-j}{k-j}+\sum_{j=1}^{k+1}N(n,j)\bi{n-j}{k-j+1}
\\=&\f1n\sum_{j=1}^{k+1}\bi nj\bi n{j-1}\bi{n-j+1}{k-j+1}
\\=&\f1n\bi nk\sum_{j=0}^{k+1}\bi nj\bi k{k+1-j}
\\=&\f1n\bi nk\bi{n+k}{k+1}=\bi nk\bi{n+k}k\f1{k+1}=[x^k]S_n(x).
\endalign$$
This proves (2.4).

The proof of Lemma 2.1 is now complete. \qed

\Remark\ 2.1. Note that the Legendre polynomial of degree $n$ is given by
$$P_n(x)=\sum_{k=0}^n\bi nk\bi{n+k}k\l(\f{x-1}2\r)^k.$$

\proclaim{Lemma 2.2} Let $n\in\Z^+$. Then
$$n(n+1)S_n(x)^2=\sum_{k=1}^n\bi{n+k}{2k}\bi{2k}k\bi{2k}{k+1}x^{k-1}(x+1)^{k+1}\tag2.5$$
and
$$\f{D_{n-1}(x)+D_{n+1}(x)}2S_n(x)=\sum_{k=0}^n\bi{n+k}{2k}\bi{2k}k^2\f{2k+1}{(k+1)^2}x^k(x+1)^{k+1}.\tag2.6$$
\endproclaim
\Remark\ 2.2. The identities (2.5) and (2.6) are (2.1) and (3.6) of the author's paper [S12a] respectively.

\proclaim{Lemma 2.3} For any $m,n\in\Z^+$ with $m\ls n$, we have the identity
$$\aligned&\sum_{k=m}^n\bi{k+m}{2m}\l(2m+1-m(m+1)\f{2k+1}{k(k+1)}\r)
\\&\qquad=\f{(n-m)(n+m+1)}{n+1}\bi{n+m}{2m}.
\endaligned\tag2.7$$
\endproclaim
\Proof. When $n=m$, both sides of (2.7) vanish.

Let $m,n\in\Z^+$ with $n\gs m$. If (2.7) holds, then
$$\align&\sum_{k=m}^{n+1}\bi{k+m}{2m}\l(2m+1-m(m+1)\f{2k+1}{k(k+1)}\r)
\\=&\f{(n-m)(n+m+1)}{n+1}\bi{n+m}{2m}
\\&+\bi{n+1+m}{2m}\l(2m+1-m(m+1)\f{2(n+1)+1}{(n+1)(n+2)}\r)
\\=&\bi{n+m+1}{2m}\l(\f{(n-m)(n-m+1)}{n+1}+2m+1-\f{m(m+1)(2n+3)}{(n+1)(n+2)}\r)
\\=&\bi{(n+1)+m}{2m}\f{(n+1-m)((n+1)+m+1)}{(n+1)+1}.
\endalign$$

In view of the above, we have proved Lemma 2.3 by induction.  \qed

Let $A$ and $B$ be integers. The Lucas sequence $u_n(A,B)\ (n=0,1,2,\ldots)$ is defined by
$u_0(A,B)=0$, $u_1(A,B)=1$, and
$$u_{n+1}(A,B)=Au_n(A,B)-Bu_{n-1}(A,B)\quad \t{for}\ n=1,2,3,\ldots.$$
It is well known that if $\Delta=A^2-4B\not=0$ then
$$u_n(A,B)=\f{\al^n-\beta^n}{\al-\beta}\quad\t{for all}\ n\in\N,$$
where $\al$ and $\beta$ are the two roots of the quadratic equation $x^2-Ax+B=0$ (so that $\al+\beta=A$ and $\alpha\beta=B$).
It is also known (see, e.g., [S10, Lemma 2.3]) that for any odd prime $p$ we have
$$u_p(A,B)\eq\l(\f{\Delta}p\r)\pmod p,\quad \t{and}\ u_{p-(\f{\Delta}p)}(A,B)\eq0\pmod{p}\ \t{if}\ p\nmid B.$$
In particular, $F_p\eq(\f 5p)\pmod p$ and $p\mid F_{p-(\f 5p)}$ for any odd prime $p$, where $F_n=u_n(1,-1)$ with $n\in\N$ is the $n$-th Fibonacci number.

\proclaim{Lemma 2.4} Let $p$ be an odd prime, and let $x$ be any integer not divisible by $p$. Then
$$W_p(x)\eq\f{4x+1}{2x}\l(\l(\f{4x+1}p\r)-1\r)\pmod p.\tag2.8$$
Moreover, if $x\eq-1/4\pmod p$ then
$$W_p(x)\eq2p\pmod{p^2},\tag2.9$$
otherwise we have
$$\aligned W_p(x)\eq &2p+\f{4x+1}{2x}\l(1-x^{p-1}+(p+1)\l(\l(\f{4x+1}p\r)-1\r)\r)
\\&-\f{4x+1}{4x^{2-(\f{4x+1}p)}}\l(2x+\l(\f{4x+1}p\r)\r)u_{p-(\f{4x+1}p)}\l(2x+1,x^2\r)
\pmod{p^2}.
\endaligned\tag2.10$$
\endproclaim
\Proof. Clearly,
$$\bi{2p-1}{p-1}=\prod_{k=1}^{p-1}\l(1+\f pk\r)\eq1+p\sum_{k=1}^{(p-1)/2}\l(\f1k+\f1{p-k}\r)\eq1\pmod{p^2}.$$
(J. Wolstenholme [W] even showed that $\bi{2p-1}{p-1}\eq1\pmod{p^3}$ if $p>3$.) Thus
$$w(p,p)=\f1p\bi{2p}{p-1}=\f2{p+1}\bi{2p-1}{p-1}\eq2(1-p)\pmod {p^2},$$
and
$$C_{p-1}=\f1p\bi{2p-2}{p-1}=\f1{2p-1}\bi{2p-1}{p-1}\eq-(2p+1)\pmod {p^2}.$$
For each $k=1,\ldots,p-1$, clearly
$$\align w(p,k)=&\f1k\bi{p-1}{k-1}\bi{p+k-1}{k-1}\f{p+k}{p+1}
\\=&\l(1+\f pk\r)\f1{p+1}\prod_{0<j<k}\l(\f{p-j}j\cdot\f{p+j}j\r)
\\\eq&\l(\f1{p+1}+\f pk\r)(-1)^{k-1}\eq(-1)^{k-1}\l(1-p+\f pk\r)\pmod{p^2}.
\endalign$$
Therefore
$$\align W_p(x)=&w(p,p)C_{p-1}x^{p-1}+\sum_{k=1}^{p-1}w(p,k)C_{k-1}x^{k-1}
\\\eq&2(1-p)C_{p-1}x^{p-1}+\sum_{k=1}^{p-1}\l(1-p+\f pk\r)C_{k-1}(-x)^{k-1}
\\\eq&(2-2p)C_{p-1}x^{p-1}+(p+1)\(\sum_{k=2}^pC_{k-1}(-x)^{k-1}+1-C_{p-1}(-x)^{p-1}\)
\\&+p\sum_{k=1}^{p-1}\l(\f1k-2\r)C_{k-1}(-x)^{k-1}
\\\eq&p+1-(1-3p)(1+2p)x^{p-1}+(p+1)\sum_{k=1}^{p-1}C_k(-x)^k
\\&-\f p2\sum_{k=1}^{p-1}\f{\bi{2k}k}k(-x)^{k-1}
\\\eq&2p+1-x^{p-1}+(p+1)\sum_{k=1}^{p-1}\f{C_k}{m^k}-\f p2m\sum_{k=1}^{p-1}\f{\bi{2k}k}{km^{k}}
\pmod {p^2},
\endalign$$
where $m$ is an integer with $m\eq-1/x\pmod{p^2}$. By [S10, Theorem 1.1], we have
$$\align \sum_{k=1}^{p-1}\f{C_k}{m^k}\eq& m^{p-1}-1-\f {m-4}2\l(\l(\f{\Delta}p\r)-1+u_{p-(\f{\Delta}p)}(m-2,1)\r)\pmod{p^2}
\\\eq&-\f{m-4}2\l(\l(\f{\Delta}p\r)-1\r)\pmod p,
\endalign$$
where
$$\Delta:=m(m-4)\eq-\f1x\l(-\f1x-4\r)=\f{4x+1}{x^2}\pmod{p^2}.$$
So (2.8) follows.

If $m\not\eq4\pmod p$, then by [S12b, Lemma 3.5] we have
$$\align \sum_{k=1}^{p-1}\f{\bi{2k}k}{km^k}\eq&\f{(-m)^{p-1}-1}p+\f{m-4}2\l(\f{\Delta}p\r)\f{u_{p-(\f{\Delta}p)}(2-m,1)}p
\\=&\f{m^{p-1}-1}p-\f{m-4}2\l(\f{\Delta}p\r)\f{u_{p-(\f{\Delta}p)}(m-2,1)}p\pmod{p}
\endalign$$
and hence from the above we deduce that
$$\align W_p(x)\eq&2p+1-x^{p-1}
\\&+(p+1)\l(m^{p-1}-1-\f {m-4}2\l(\l(\f{\Delta}p\r)-1+u_{p-(\f{\Delta}p)}(m-2,1)\r)\r)
\\&-\f p2m\l(\f{m^{p-1}-1}p-\f{m-4}2\l(\f{\Delta}p\r)\f{u_{p-(\f{\Delta}p)}(m-2,1)}p\r)
\\\eq&2p+1-x^{p-1}+\l(1-\f m2\r)(m^{p-1}-1)-(p+1)\f{m-4}2\l(\l(\f{4x+1}p\r)-1\r)
\\&-\f{m-4}2\l(1-\f m2\l(\f{4x+1}p\r)\r)u_{p-(\f{4x+1}p)}(m-2,1)
\\\eq&2p+1-x^{p-1}+\f{2x+1}{2x}(1-x^{p-1})+(p+1)\f{4x+1}{2x}\l(\l(\f{4x+1}p\r)-1\r)
\\&+\f{4x+1}{2x}\l(1+\f1{2x}\l(\f{4x+1}p\r)\r)u_{p-(\f{4x+1}p)}(m-2,1)\pmod{p^2}
\endalign$$
which is equivalent to (2.10) since
$$(-x)^{k-1}u_k(m-2,1)=u_k(-x(m-2),(-x)^2)\eq u_k(2x+1,x^2)\pmod{p^2}$$
for all $k\in\N$.

When $m\eq4\pmod p$ (i.e., $x\eq-1/4\pmod p$), we have
$$\sum_{k=1}^{p-1}\f{\bi{2k}k}{km^k}\eq\sum_{k=1}^{p-1}\f{\bi{2k}k}{k4^k}\eq2q_p(2)\pmod p$$
by [ST, (1.12)], and hence
$$\align W_p(x)\eq& 2p+1-x^{p-1}+(p+1)\l(m^{p-1}-1+\f{m-4}2\r)-\f p2m2q_p(2)
\\\eq&2p+1-x^{p-1}+m^{p-1}-1+\f{m-4}2-4(2^{p-1}-1)
\\\eq&2p+1-x^{p-1}+1-x^{p-1}-\f{4x+1}{2x}-2(2^{2(p-1)}-1)
\\\eq&2\l(p+(4x+1)+(1-x^{p-1})+(1-4^{p-1})\r)
\\\eq&2\l(p+(4x+1)+1-(4x+1-1)^{p-1}\r)
\\\eq&2\l(p+(4x+1)+(p-1)(4x+1)\r)\eq2p\pmod{p^2}.
\endalign$$
Therefore (2.9) is valid. \qed

\medskip
\noindent{\it Proof of Theorem 1.1}. (i) Fix $n\in\Z^+$.
For each $k\in\Z^+$, by (2.3) and Lemma 2.2 we have
$$\align D_{k-1}(x)\f{S_k(x)}{x+1}=&\f{D_{k-1}(x)+D_{k+1}(x)}2\cdot\f{S_k(x)}{x+1}-(2k+1)\f x{x+1}S_k(x)^2
\\=&\sum_{j=0}^k\bi{k+j}{2j}\bi{2j}j^2\f{2j+1}{(j+1)^2}(x(x+1))^j
\\&-\f{2k+1}{k(k+1)}\sum_{j=0}^k\bi{k+j}{2j}\bi{2j}j\bi{2j}{j+1}(x(x+1))^j
\\=&\sum_{j=0}^k\bi{k+j}{2j}\bi{2j}j^2\l(\f{2j+1}{(j+1)^2}-\f{2k+1}{k(k+1)}\cdot\f j{j+1}\r)(x(x+1))^j
\\=&\sum_{j=0}^k\bi{k+j}{2j}C_j^2\l(2j+1-j(j+1)\f{2k+1}{k(k+1)}\r)(x(x+1))^j.
\endalign$$
Combining this with (2.4) we obtain
$$D_{k-1}(x)s_k(x)=\sum_{j=0}^k\bi{k+j}{2j}C_j^2\l(2j+1-j(j+1)\f{2k+1}{k(k+1)}\r)(x(x+1))^j\tag2.11$$
for any $k\in\Z^+$. Therefore,
$$\align \sum_{k=1}^{n}D_{k-1}(x)s_{k}(x)
=&\sum_{k=1}^n\sum_{j=0}^k\bi{k+j}{2j}C_j^2\l(2j+1-j(j+1)\f{2k+1}{k(k+1)}\r)(x(x+1))^j
\\=&n+\sum_{j=1}^nC_j^2(x(x+1))^j\sum_{k=j}^n\bi{k+j}{2j}\l(2j+1-j(j+1)\f{2k+1}{k(k+1)}\r)
\\=&\sum_{j=0}^{n-1}C_j^2(x(x+1))^j\f{(n-j)(n+j+1)}{n+1}\bi{n+j}{2j}
\endalign$$
with the help of Lemma 2.3. It follows that
$$\sum_{k=0}^{n-1}D_k(x)s_{k+1}(x)=\sum_{j=0}^{n-1}C_j(x(x+1))^j\f n{j+1}\bi{n-1}j\bi{n+j+1}j
=nW_n(x(x+1)).$$
For any $k\in\Z^+$, we have
$$w(n,k)=\f1n\bi nk\bi{n+k}{k-1}=\f1{n+1}\bi{n-1}{k-1}\bi{n+k}k$$
and hence $w(n,k)=(n+1)w(n,k)-nw(n,k)\in\Z$.
So $W_n(x)\in\Z[x]$. This proves part (i) of Theorem 1.1.

(ii) Let $x$ be any $p$-adic integer, and let $y$ be an integer with $y\eq x(x+1)\pmod{p^2}$.
By (1.7) we have
$$\sum_{k=0}^{p-1}D_k(x)s_{k+1}(x)=pW_p(x(x+1))\eq pW_p(y)\pmod{p^3}.\tag2.12$$

If $y\eq0\pmod p$, then
$$W_p(y)\eq w(p,1)+w(p,2)y\eq1-x(x+1)\pmod{p^2}.$$
If $x\eq-1/2\pmod p$ (i.e., $y\eq-1/4\pmod p$), then
$W_p(y)\eq2p\pmod{p^2}$ by Lemma 2.4.
Thus (1.10) holds for $x\eq0,-1,-1/2\pmod p$.

Now assume that $x\not\eq0,-1,-1/2\pmod{p}$, i.e., $y\not\eq0,-1/4\pmod p$. Then
$$\l(\f{4y+1}p\r)=\l(\f{(2x+1)^2}p\r)=1$$
and
$$\align u_{p-1}(2y+1,y^2)\eq&u_{p-1}(x^2+(x+1)^2,x^2(x+1)^2)
\\=&\f{((x+1)^2)^{p-1}-(x^2)^{p-1}}{(x+1)^2-x^2}
\\=&\f{(x+1)^{p-1}+x^{p-1}}{2x+1}\l((x+1)^{p-1}-x^{p-1}\r)
\\\eq&\f2{2x+1}\l((x+1)^{p-1}-x^{p-1}\r)\pmod{p^2}.
\endalign$$
Combining this with Lemma 2.4, we obtain
$$\align W_p(y)\eq&2p+\f{(2x+1)^2}{2x(x+1)}\l(1-x^{p-1}(x+1)^{p-1}\r)
\\&-\f{(2x+1)^2}{4x(x+1)}(2x(x+1)+1)\f2{2x+1}\l((x+1)^{p-1}-x^{p-1}\r)
\\\eq&2p+\f{(2x+1)^2}{2x(x+1)}\l(1-x^{p-1}+1-(x+1)^{p-1}\r)
\\&-\f{2x+1}{2x(x+1)}(2x^2+2x+1)\l((x+1)^{p-1}-x^{p-1}\r)
\\=&2p+\f{2x+1}{x(x+1)}\l(x^2(x^{p-1}-1)-(x+1)^2\l((x+1)^{p-1}-1\r)\r)
\pmod {p^2}.
\endalign$$
Therefore (1.10) holds in light of (2.12).

So far we have completed the proof of Theorem 1.1. \qed

\medskip
\noindent{\it Proof of Corollary 1.2}.
It is known (cf. [G, (3.133)]) that
$$D_k\l(-\f12\r)=P_k(0)=\cases (-1)^{k/2}\bi k{k/2}/2^k&\t{if}\ 2\mid k,
\\0&\t{otherwise}.\endcases$$
By (2.4) and [S11a, Lemma 4.3],
$$s_{k+1}\l(-\f12\r)=2S_{k+1}\l(-\f12\r)=\cases(-1)^{k/2}C_{k/2}/2^k&\t{if}\ 2\mid k,\\0&\t{otherwise}.\endcases$$
Therefore
$$\align&\sum_{k=0}^{n-1}D_k\l(-\f12\r)s_{k+1}\l(-\f12\r)
\\=&\sum\Sb 0\ls k<n\\2\mid k\endSb\l(\f{(-1)^{k/2}}{2^k}\r)^2\bi{k}{k/2}C_{k/2}=\sum_{j=0}^{\lfloor(n-1)/2\rfloor}\f{\bi{2j}j^2}{(j+1)16^j}.
\endalign$$
By induction,
$$\sum_{j=0}^m\f{\bi{2j}j^2}{(j+1)16^j}=\f{(2m+1)^2}{(m+1)16^m}\bi{2m}m^2\quad\t{for all}\ m\in\N.$$
So we have
$$\align&\sum_{k=0}^{n-1}D_k\l(-\f12\r)s_{k+1}\l(-\f12\r)
\\=&\f{(2\lfloor(n-1)/2\rfloor+1)^2}{\lfloor(n+1)/2\rfloor16^{\lfloor(n-1)/2\rfloor}}\bi{2\lfloor(n-1)/2\rfloor}{\lfloor(n-1)/2\rfloor}^2
=\f{\lfloor(n+1)/2\rfloor}{4^{n-1}}\bi n{\lfloor n/2\rfloor}^2.
\endalign$$
On the other hand, by applying (1.7) with $x=-1/2$ we obtain
$$\sum_{k=0}^{n-1}D_k\l(-\f12\r)s_{k+1}\l(-\f12\r)=nW_n\l(-\f14\r)=\sum_{k=1}^n\bi nk\bi{n+k}{k-1}\f{C_{k-1}}{(-4)^{k-1}}.$$
Combining these we get the desired identity (1.13). This concludes the proof. \qed

\proclaim{Lemma 2.5} Let $p$ be any odd prime.
For each $j=0,\ldots,p$, we have
$$C_j^2u_j\eq\cases2\pmod{p}&\t{if}\ j=(p-1)/2,\\0\pmod p&\t{otherwise},\endcases\tag2.13$$
where
$$u_j:=\sum_{j<k\ls p}(k-1)\bi{k+j}{2j}\l(2j+1-j(j+1)\f{2k+1}{k(k+1)}\r).\tag2.14$$
\endproclaim
\Proof. Clearly, $u_p=0$ and
$$u_{p-1}=(p-1)\bi{p+(p-1)}{2(p-1)}\l(2(p-1)+1-(p-1)p\f{2p+1}{p(p+1)}\r)\eq0\pmod p.$$
If $(p-1)/2<j<p-1$, then $C_j=(2j)!/(j!(j+1)!)\eq0\pmod p$ and hence
$C_j^2u_j\eq0\pmod p$ since $pu_j$ is a $p$-adic integer.

Note that
$$C_{(p-1)/2}=\f2{p+1}\bi{p-1}{(p-1)/2}\eq2(-1)^{(p-1)/2}\pmod p.$$
As
$$\bi{k+(p-1)/2}{p-1}\eq0\pmod p\quad\t{for all}\ k=\f{p+1}2,\ldots,p,$$
we have
$$\align u_{(p-1)/2}\eq&\sum_{k=p-1}^p(k-1)\bi{k+(p-1)/2}{p-1}\l(2\cdot\f{p-1}2+1-\f{p-1}2\cdot\f{p+1}2\cdot\f{2k+1}{k(k+1)}\r)
\\\eq&\f14\sum_{k=p-1}^p(k-1)\bi{k+(p-1)/2}{p-1}\f{2k+1}{k(k+1)}
\\=&\f{p-2}4\bi{p-1+(p-1)/2}{(p-1)/2}\f{2(p-1)+1}{(p-1)p}+\f{p-1}4\bi{p+(p-1)/2}{(p+1)/2}\f{2p+1}{p(p+1)}
\\=&\f{p-2}4\cdot\f{2p-1}{p-1}\cdot\f{\prod_{1<r\ls(p-1)/2}(p-1+r)}{((p-1)/2)!}
\\&+\f{p-1}4\cdot\f{2p+1}{p+1}\cdot\f{\prod_{r=1}^{(p-1)/2}(p+r)}{((p+1)/2)!}
\\\eq&-\f12\cdot\f1{(p-1)/2}-\f14\cdot\f1{(p+1)/2}\eq\f12\pmod p.
\endalign$$

Obviously,
$$u_0=\sum_{k=1}^p(k-1)=\f{p(p-1)}2\eq0\pmod p.$$
Applying the Zeilberger algorithm via {\tt Mathematica 9}, we find that
$$(j+2)u_j+2(2j+1)u_{j+1}=\f{f(p,j)\bi{p+j}{2j}}{2(j+1)(j+2)(2j+3)}$$
for all $j=0,\ldots,p-1$, where
$$f(p,j):=(p-j)(p+j+1)\l((2j+3)^2p^2-(2j^2+8j+7)p-(j+1)(j+2)\r).$$
This implies that for $0\ls j<(p-3)/2$ we have
$$u_j\eq0\pmod p \Longrightarrow u_{j+1}\eq0\pmod p.$$
Thus $u_j\eq0\pmod p$ for all $j=0,\ldots,(p-3)/2$.

Combining the above, we immediately obtain the desired (2.13). \qed

\medskip
\noindent{\it Proof of Theorem 1.2}. In view of (2.11),
$$\align&\sum_{k=1}^p(k-1)D_{k-1}(x)s_k(x)
\\=&\sum_{k=1}^p(k-1)\sum_{j=0}^k\bi{k+j}{2j}C_j^2\l(2j+1-j(j+1)\f{2k+1}{k(k+1)}\r)(x(x+1))^j
\\=&\sum_{k=1}^p(k-1)+\sum_{j=1}^p(x(x+1))^jC_j^2u_j=\f{p(p-1)}2+\sum_{j=1}^pC_j^2u_j(x(x+1))^j,
\endalign$$
where $u_j$ is given by (2.14). Thus, by applying Lemma 2.5 we find that
$$\f1p\(\sum_{k=0}^{p-1}kD_k(x)s_{k+1}(x)-2(x(x+1))^{(p-1)/2}\)\in\Z_p[x(x+1)],\tag2.15$$
where $\Z_p$ denotes the ring of $p$-adic integers.
Therefore (1.14) holds. (1.14) with $x=1$ gives (1.15). This concludes the proof. \qed

\heading{3. Applications to central trinomial coefficients and Motzkin numbers}\endheading

Let $n\in\N$ and $b,c\in\Z$. The $n$-th generalized central trinomial coefficient $T_n(b,c)$ is defined to be $[x^n](x^2+bx+c)^n$, the
coefficient of $x^n$ in the expansion of $(x^2+bx+c)^n$. It is easy to see that
$$T_n(b,c)=\sum_{k=0}^{\lfloor n/2\rfloor}\bi n{2k}\bi{2k}kb^{n-2k}c^k.\tag3.1$$
Note that $T_n(1,1)$ is the central trinomial coefficient $T_n$ and $T_n(2,1)$ is the central binomial coefficient $\bi{2n}n$.
Also, $T_n(3,2)$ coincides with the central Delannoy number $D_n$.
Sun [S14a] also defined the generalized Motzkin number $M_n(b,c)$ by
$$M_n(b,c)=\sum_{k=0}^{\lfloor n/2\rfloor}\bi n{2k}C_kb^{n-2k}c^k.\tag3.2$$
Note that $M_n(1,1)$ is the usual Motzkin number $M_n$ (whose combinatorial interpretations can be found in [St99, Ex.\,6.38])
 and $M_n(2,1)$ is the Catalan number $C_{n+1}$.
Also, $M_n(3,2)$ coincides with the little Schr\"oder number $s_{n+1}$.
The author [S14a, S14b] deduced some congruences involving $T_n(b,c)$ and $M_n(b,c)$, and proposed in [S14b]
some conjectural series for $1/\pi$ involving $T_n(b,c)$ such as
$$\align
\sum_{k=0}^\infty\f{66k+17}{(2^{11}3^3)^{k}}T_k(10,11^2)^3=&\f{540\sqrt2}{11\pi},
\\\sum_{k=0}^\infty\f{126k+31}{(-80)^{3k}}T_k(22,21^2)^3=&\f{880\sqrt5}{21\pi},
\\\sum_{k=0}^\infty\f{3990k+1147}{(-288)^{3k}}T_k(62,95^2)^3=&\f{432}{95\pi}(195\sqrt{14}+94\sqrt2).
\endalign$$

Now we point out that $T_n(b,c)$ and $M_n(b,c)$ are actually related to the polynomials $D_n(x)$ and $s_{n+1}(x)$.

\proclaim{Lemma 3.1} Let $b,c\in\Z$ with $d=b^2-4c\not=0$. For any $n\in\N$ we have
$$T_n(b,c)=(\sqrt d)^nD_n\l(\f{b/\sqrt d-1}2\r)\tag3.3$$
and
$$M_n(b,c)=(\sqrt d)^ns_{n+1}\l(\f{b/\sqrt d-1}2\r).\tag3.4$$
\endproclaim
\Proof. In view of (3.1) and (3.2),
$$\f{T_n(b,c)}{(\sqrt d)^n}=\sum_{k=0}^{\lfloor n/2\rfloor}\bi n{2k}\bi{2k}k\l(\f b{\sqrt d}\r)^{n-2k}\l(\f cd\r)^k$$
and
$$\f{M_n(b,c)}{(\sqrt d)^n}=\sum_{k=0}^{\lfloor n/2\rfloor}\bi n{2k}C_k\l(\f b{\sqrt d}\r)^{n-2k}\l(\f cd\r)^k.$$
Note that
$$\l(\f{b}{\sqrt d}\r)^2-4\f cd=1. $$
So, it suffices to show the polynomial identities
$$\sum_{k=0}^{\lfloor n/2\rfloor}\bi n{2k}\bi{2k}k(2x+1)^{n-2k}(x(x+1))^k= D_n(x)\tag3.5$$
and
$$\sum_{k=0}^{\lfloor n/2\rfloor}\bi n{2k}C_k(2x+1)^{n-2k}(x(x+1))^k= s_{n+1}(x).\tag3.6$$

It is easy to verify (3.5) and (3.6) for $n=0,1,2$.
Let $u_n(x)$ denote the left-hand side or the right-hand side of $(3.5)$. By the Zeilberger algorithm (cf. [PWZ, pp.\,101-119]), we have the recurrence
$$(n+2)u_{n+2}(x)=(2x+1)(2n+3)u_{n+1}(x)-(n+1)u_n(x)\ \ \t{for}\ n=0,1,2,\ldots.$$
Thus (3.5) is valid by induction.
Let $v_n(x)$ denote the left-hand side or the right-hand side of $(3.6)$. By the Zeilberger algorithm, we have the recurrence
$$(n+4)v_{n+2}(x)=(2x+1)(2n+5)v_{n+1}(x)-(n+1)v_n(x)\ \ \t{for}\ n=0,1,2,\ldots.$$
So (3.6) also holds by induction.

The proof of Lemma 3.1 is now complete. \qed

With the help of Theorem 1.1, we are able to confirm Conjecture 5.5 of the author [S14a] by proving the following result.

\proclaim{Theorem 3.1} Let $b,c\in\Z$ and $d=b^2-4c$.

{\rm (i) For any $n\in\Z^+$, we have
$$\f1n\sum_{k=0}^{n-1}T_k(b,c)M_k(b,c)d^{n-1-k}=\sum_{k=1}^nw(n,k)C_{k-1}c^{k-1}d^{n-k}\in\Z.\tag3.7$$
Moreover, for any odd prime $p$ not dividing $cd$, we have
$$\sum_{k=0}^{p-1}\f{T_k(b,c)M_k(b,c)}{d^k}\eq\f{pb^2}{2c}\l(\l(\f dp\r)-1\r)\pmod{p^2},\tag3.8$$
and furthermore
$$\aligned&\sum_{k=0}^{p-1}\f{T_k(b,c)M_k(b,c)}{d^k}
\\\eq&\f{pb^2}{2c}\l(\l(\f dp\r)-1\r)+\f{p^2}{2c}\l(b^2\l(q_p(d)-q_p(c)+\l(\f dp\r)\r)-d\r)
\\&-\f{pb^2}{4c^{2-(\f dp)}}\l(2c+d\l(\f dp\r)\r)u_{p-(\f dp)}(b^2-2c,c^2)\pmod{p^3}.
\endaligned\tag3.9$$

{\rm (ii)}  For any odd prime $p$ not dividing $d$, we have
$$\sum_{k=0}^{p-1}\f{kT_k(b,c)M_k(b,c)}{d^k}\eq2\l(\f {cd}p\r)\pmod p.\tag3.10$$
\endproclaim
\Proof. (i) Let's first prove (3.7) for any $n\in\Z^+$.

We first consider the case $d=0$, i.e., $c=b^2/4$.
In this case, for any $k\in\N$ we have
$$T_k(b,c)=\l(\f b2\r)^kT_k(2,1)=\l(\f b2\r)^k\bi{2k}k$$
and
$$M_k(b,c)=\l(\f b2\r)^kM_k(2,1)=\l(\f b2\r)^kC_{k+1}.$$
Thus
$$\align&\f1n\sum_{k=0}^{n-1}T_k(b,c)M_k(b,c)d^{n-1-k}
\\=&\f{T_{n-1}(b,c)M_{n-1}(b,c)}n=\f1n\l(\f{b^2}4\r)^{n-1}\bi{2(n-1)}{n-1}C_n
\\=&c^{n-1}C_{n-1}C_n=w(n,n)C_{n-1}c^{n-1}
\endalign$$
and hence (3.7) is valid.

Now assume that $d\not=0$. By Lemma 3.1 and (1.7), we have
$$\align&\f1n\sum_{k=0}^{n-1}\f{T_k(b,c)M_k(b,c)}{d^k}
\\=&\f1n\sum_{k=0}^{n-1}D_k\l(\f{b/\sqrt d-1}2\r)s_{k+1}\l(\f{b/\sqrt d-1}2\r)
\\=&W_n\l(\f{b/\sqrt d-1}2\cdot\f{b/\sqrt d+1}2\r)=W_n\l(\f{b^2/d-1}4\r)=W_n\l(\f cd\r)
\endalign$$
and hence (3.7) holds in view of (1.8).

Below we suppose that $p$ is an odd prime not dividing $cd$.
From the above, we have
$$\sum_{k=0}^{p-1}\f{T_k(b,c)M_k(b,c)}{d^k}=pW_p\l(\f cd\r)\eq pW_p(x)\pmod{p^3},\tag3.11$$
where $x$ is an integer with $x\eq c/d\pmod{p^2}$. As $p\nmid d$ and $d(4x+1)\eq 4c+d=b^2\pmod{p^2}$, we have
$$\l(\f{4x+1}p\r)=\l(\f{d^2(4x+1)}p\r)=\l(\f{b^2d}p\r).$$
In view of Lemma 2.4,
$$W_p\l(x\r)\eq\f{4c/d+1}{2c/d}\l(\l(\f{4x+1}p\r)-1\r)\eq\f{b^2}{2c}\l(\l(\f{d}p\r)-1\r)\pmod p.$$
Combining this with (3.11) we immediately obtain (3.8).

Now we show (3.9). If $x\eq-1/4\pmod p$ (i.e., $p\mid b$), then by (2.9) we have
$$W_p\l(x\r)\eq2p\eq\f p{2c}(4c-b^2)\pmod {p^2}$$
and hence (3.9) holds by (3.11). Below we assume $p\nmid b$. Then
$$\l(\f{(2x+1)^2-4x^2}p\r)=\l(\f{4x+1}p\r)=\l(\f{b^2d}p\r)=\l(\f dp\r),$$
$p\mid u_{p-(\f dp)}(2x+1,x^2)$ and
$$\align &d^{p-1-(\f dp)}u_{p-(\f dp)}\l(2x+1,x^2\r)
\\=&u_{p-(\f dp)}(d(2x+1),d^2x^2)
\\\eq&u_{p-(\f dp)}(2c+d,c^2)=u_{p-(\f dp)}(b^2-2c,c^2)\pmod{p^2}.
\endalign$$
So, applying (2.10) we get
$$\align W_p\l(x\r)\eq&2p+\f{4c/d+1}{2c/d}\l(1-\l(\f cd\r)^{p-1}+(p+1)\l(\l(\f dp\r)-1\r)\r)
\\&-\f{4c/d+1}{4(c/d)^{2-(\f dp)}}\l(2\f cd+\l(\f dp\r)\r)d^{(\f dp)}u_{p-(\f dp)}(b^2-2c,c^2)
\\\eq&2p+\f{b^2}{2c}\l(d^{p-1}-1-(c^{p-1}-1)+(p+1)\l(\l(\f dp\r)-1\r)\r)
\\&-\f{b^2}{4c^{2-(\f dp)}}\l(2c+d\l(\f dp\r)\r)u_{p-(\f dp)}(b^2-2c,c^2)\pmod{p^2}.
\endalign$$
This, together with (3.11), yields the desired (3.9).

(ii) Fix an odd prime $p$ not dividing $d$. Let $x=b/\sqrt d-1$. Then
$$x(x+1)=\f{b/\sqrt d-1}2\cdot\f{b/\sqrt d+1}2=\f{b^2/d-1}4=\f cd$$
is a $p$-adic integer. Thus, with the help of (2.15), we have
$$\sum_{k=0}^{p-1}kD_k(x)s_{k+1}(x)\eq 2\l(\f cd\r)^{(p-1)/2}\eq2\l(\f{cd}p\r)\pmod p.$$
Combining this with Lemma 3.1, we immediately obtain (3.10).

In view of the above, we have proved Theorem 3.1. \qed
\medskip

Let $\omega$ denote the primitive cubic root $(-1+\sqrt{-3})/2$ of unity.
Then $\omega+\bar\omega=-1$ and $\omega\bar\omega=1$. So,
$$u_n(-1,1)=\f{\omega^n-\bar\omega^n}{\omega-\bar\omega}=0\ \ \t{and}\ \ u_n(3,9)=\f{(-3\omega)^n-(-3\bar\omega)^n}{(-3\omega)-(-3\bar\omega)}=0$$
for any $n\in\N$ with $3\mid n$.
In view of this, Theorem 3.1 in the cases $b=c\in\{1,3\}$ yields the following consequence.

\proclaim{Corollary 3.1} For any positive integer $n$, we have
$$\f1n\sum_{k=0}^{n-1}T_kM_k(-3)^{n-1-k}=\sum_{k=1}^nw(n,k)C_{k-1}(-3)^{n-k}\in\Z\tag3.12$$
and
$$\f1n\sum_{k=0}^{n-1}\f{T_k(3,3)M_k(3,3)}{(-3)^k}=\sum_{k=1}^n(-1)^{k-1}w(n,k)C_{k-1}\in\Z.\tag3.13$$
Moreover, for any prime $p>3$ we have
$$\sum_{k=0}^{p-1}\f{T_kM_k}{(-3)^k}\eq\f p2\l(\l(\f{p}3\r)-1\r)+\f{p^2}2\l(q_p(3)+\l(\f p3\r)+3\r)\pmod{p^3},\tag3.14$$
$$\sum_{k=0}^{p-1}\f{T_k(3,3)M_k(3,3)}{(-3)^k}\eq\f{3p}2\l(\l(\f p3\r)-1\r)+\f{p^2}2\l(3\l(\f p3\r)+1\r)\pmod{p^3},\tag3.15$$
$$\sum_{k=0}^{p-1}\f{kT_kM_k}{(-3)^k}\eq2\l(\f p3\r)\ (\mo\ p) \ \t{and}\ \sum_{k=0}^{p-1}\f{kT_k(3,3)M_k(3,3)}{(-3)^k}\eq2\l(\f {-1}p\r)\ (\mo\ p).\tag3.16$$
\endproclaim
\Remark\ 3.1. Let $p>3$ be a prime. In the case $p\eq2\pmod3$, the author (cf. [S14a, Conjecture 5.6]) even conjectured that
$$\sum_{k=0}^{p-1}\f{T_k(3,3)M_k(3,3)}{(-3)^k}\eq p^3-p^2-3p\pmod{p^4}$$
which is stronger than (3.15).
The author's conjectural supercongruences (cf. [S14a, Conjecture 1.1(ii)])
$$\sum_{k=0}^{p-1}M_k^2\eq(2-6p)\l(\f p3\r)\pmod{p^2},\ \ \sum_{k=0}^{p-1}kM_k^2\eq(9p-1)\l(\f p3\r)\pmod{p^2},$$
and
$$\sum_{k=0}^{p-1}T_kM_k\eq\f 43\l(\f p3\r)+\f p6\l(1-9\l(\f p3\r)\r)\pmod{p^2}$$
remain open. We also observe that
$$\sum_{k=0}^{p-1}kT_kM_k\eq\l(\f{-1}p\r)-\f 53\l(\f p3\r)\pmod p.$$
\medskip

The Lucas numbers $L_0,L_1,L_2,\ldots$ are given by
$$L_0=2,\ L_1=1,\ \t{and}\ L_{n+1}=L_n+L_{n-1}\ \t{for}\ n=1,2,3,\ldots.$$
It is easy to see that $L_n=2F_{n+1}-F_n=2F_{n-1}+F_n$ for all $n\in\Z^+$.
Thus, for any odd prime $p\not=5$ we have
$$L_{p-(\f p5)}=2F_p-\l(\f p5\r)F_{p-(\f p5)}\eq2\l(\f p5\r)\pmod p$$
and hence
$$u_{p-(\f p5)}(3,1)=\f{(\al^2)^{p-(\f p5)}-(\beta^2)^{p-(\f p5)}}{\al^2-\beta^2}=F_{p-(\f p5)}L_{p-(\f p5)}\eq2\l(\f p5\r)F_{p-(\f p5)}\ (\mo\ {p^2}),$$
where $\al=(1+\sqrt 5)/2$ and $\beta=(1-\sqrt5)/2$.
Note also that
$$u_n(3\times 5,5^2)=5^{n-1}u_n(3,1)=5^{n-1}F_nL_n\quad\t{for any}\ n\in\N.$$
Thus Theorem 3.1 with $(b,c)=(1,-1),(5,5)$ leads to the following corollary.

\proclaim{Corollary 3.2} For any $n\in\Z^+$, we have
$$\f1n\sum_{k=0}^{n-1}T_k(1,-1)M_k(1,-1)5^{n-1-k}=\sum_{k=1}^n(-1)^{k-1}w(n,k)C_{k-1}5^{n-k}\in\Z\tag3.17$$
and
$$\f1n\sum_{k=0}^{n-1}\f{T_k(5,5)M_k(5,5)}{5^k}=\sum_{k=1}^nw(n,k)C_{k-1}\in\Z.\tag3.18$$
Also, for any prime $p\not=2,5$ we have the congruences
$$\aligned \sum_{k=0}^{p-1}\f{T_k(1,-1)M_k(1,-1)}{5^k}
\eq&\f p2\l(1-\l(\f p5\r)\r)+\f{p^2}2\l(5-\l(\f p5\r)-q_p(5)\r)
\\&+\f p2\l(5-2\l(\f p5\r)\r)F_{p-(\f p5)}\pmod{p^3},
\endaligned\tag3.19$$
$$\aligned \sum_{k=0}^{p-1}\f{T_k(5,5)M_k(5,5)}{5^k}
\eq&\f {5p}2\l(\l(\f p5\r)-1\r)+\f{p^2}2\l(5\l(\f p5\r)-1\r)
\\&-\f {5p}2\l(1+2\l(\f p5\r)\r)F_{p-(\f p5)}\pmod{p^3},
\endaligned\tag3.20$$
and
$$\l(\f{-5}p\r)\sum_{k=0}^{p-1}\f{kT_k(1,-1)M_k(1,-1)}{5^k}\eq\sum_{k=0}^{p-1}\f{kT_k(5,5)M_k(5,5)}{5^k}\eq2\pmod p.\tag3.21$$
\endproclaim

\heading{4. Two related conjectures}\endheading

In view of (1.8) and (1.9), we can easily see that
$$\align&W_1(x)=1,\ W_2(x)=2x+1,\ W_3(x)=10x^2+5x+1,
\\& W_4(x)=70x^3+42x^2+9x+1.\endalign$$
Applying the Zeilberger algorithm (cf. [PWZ, pp.\,101-119]) via {\tt Mathematica 9}, we obtain the following third-order recurrence with $n\in\Z^+$:
$$\aligned &(n+3)^2(n+4)(2n+3)W_{n+3}(x)
\\=&(n+3)(2n+5)(4x(2n+3)^2+3n^2+11n+10)W_{n+2}(x)
\\&-(n+1)(2n+3)(4x(2n+5)^2+3n^2+13n+14)W_{n+1}(x)
\\&+n(n+1)^2(2n+5)W_n(x).
\endaligned\tag4.1$$
(This is a verified result, not a conjecture.)

For any $n\in\Z^+$, we clearly have $w(n,n)=C_n$.
For the polynomial
$$w_n(x):=\sum_{k=1}^nw(n,k)x^{k-1},\tag4.2$$
we have the relation
$$w_n(-1-x)=(-1)^{n-1}w_n(x)\tag4.3$$
since
$$\sum_{k=m}^n(-1)^{n-k}\bi{k-1}{m-1}w(n,k)=w(n,m)\ \quad\t{for all}\ m=1,\ldots,n,\tag4.4$$
which can be deduced with the help of the Chu-Vandermonde identity in the following way:
$$\align &\sum_{k=m}^n(-1)^{n-k}\bi{k-1}{m-1}w(n,k)
\\=&\bi{n-1}{m-1}\sum_{k=m}^n\f{(-1)^{n-k}}k\bi{n-m}{n-k}\bi{-n-2}{k-1}(-1)^{k-1}
\\=&\f{(-1)^n}{n+1}\bi{n-1}{m-1}\sum_{k=m}^n\bi{n-m}{n-k}\bi{-n-1}k
\\=&\f{(-1)^n}{n+1}\bi{n-1}{m-1}\bi{-m-1}n=w(n,m).
\endalign$$
Via the Zeilberger algorithm we obtain the recurrence
$$(n+3)w_{n+2}(x)=(2x+1)(2n+3)w_{n+1}(x)-nw_n(x)\quad\t{for}\ n=1,2,3,\ldots\tag4.5$$
As $w_2(x)=2x+1$, this recurrence implies that $w_{2n}(-1/2)=0$
and hence $w_{2n}(x)/(2x+1)\in\Z[x]$ for all $n\in\Z^+$. We also note that
$$\sum_{n=1}^\infty w_n(x)y^{n}=\f{1-y-2xy-\sqrt{(y-1)^2-4xy}}{2x(x+1)y},\tag4.6$$
while
$$\sum_{n=0}^\infty S_ny^n=\f{1-y-\sqrt{y^2-6y+1}}{2y}.$$
\medskip

Now we pose two conjectures for further research.

\proclaim{Conjecture 4.1} For any integer $n>1$, all the polynomials
$$w_{2n-1}(x),\ \ \f{w_{2n}(x)}{2x+1}\ \ \t{and}\ \ W_n(x)$$
are irreducible over the field of rational numbers.
\endproclaim

\proclaim{Conjecture 4.2} {\rm (i)} For any $n\in\Z^+$, we have
$$f_n(x):=\f1n\sum_{k=0}^{n-1}D_k(x)R_k(x)\in\Z[x],\tag4.7$$
where
$$R_k(x):=\sum_{l=0}^k\bi kl\bi{k+l}l\f{x^l}{2l-1}=\sum_{l=0}^k\bi{k+l}{2l}\bi{2l}l\f{x^l}{2l-1}.\tag4.8$$
Also, $f_2(x),f_3(x),\ldots$ are all irreducible over the field of rational numbers,
and
$$f_n(1)=\f1n\sum_{k=0}^{n-1}D_kR_k\eq(-1)^n\pmod{32}$$
  for each $n\in\Z^+$, where $R_k=R_k(1)$.

{\rm (ii)} Let $p$ be any odd prime. Then
$$\sum_{k=0}^{p-1}D_kR_k\eq\cases-p+8p^2q_p(2)-2p^3E_{p-3}\pmod{p^4}&\t{if}\ p\eq1\pmod 4,
\\-5p\pmod{p^3}&\t{if}\ p\eq3\pmod4.\endcases\tag4.9$$
Also,
$$\align\sum_{k=1}^{p-1}\f{D_kR_k}k\eq&\l(4-\l(\f{-1}p\r)\r)q_p(2)\pmod p,\tag4.10
\\\sum_{k=1}^{p-1}kD_kR_k\eq&\f12+\f32p\l(1-2\l(\f{-1}p\r)\r)\pmod{p^2},\tag4.11
\endalign$$
and
$$\sum_{k=1}^{p-1}kD_k(x)R_k(x)\eq\f{x^{p-1}}2\pmod p.$$
\endproclaim
\medskip

\Ack. The author would like to thank the referee for helpful comments.

\enddocument